\documentclass{amsart}
\usepackage{amsfonts,amssymb,amsmath,a4wide}
\usepackage[british]{babel}

\newtheorem{teo}{Theorem}[section]
\newtheorem{lema}{Lemma}[section]
\newtheorem{pro}{Proposition}[section]
\newtheorem{defi}{Definition}[section]
\newtheorem{rema}{Remark}[section]
\newtheorem{coro}{Corollary}[section]
\numberwithin{equation}{section}

\def \z{\zeta}
\def \S {\mathcal{S}}
\def \E {\mathcal{E}}
\def \A {\mathcal{A}}
\def \W{\mathcal{W}}

\begin{document}\larger[2]
\title{Four dimensional conformal $C$-spaces}
\author[A.R.Gover]{A. Rod Gover}
\address[A.R.Gover]{Department of Mathematics, University of 
Auckland, Private Bag 92019, Auckland, New Zealand}
\email{gover@math.auckland.ac.nz}
\author[P.-A. Nagy]{Paul-Andi Nagy}
\address[P.-A. Nagy]{Department of Mathematics, University of 
Auckland, Private Bag 92019, Auckland, New Zealand}
\email{nagy@math.auckland.ac.nz}
\date{\today}
\begin{abstract} We investigate the structure of conformal $C$-spaces,
  a class of Riemmanian manifolds which naturally arises as a
  conformal generalisation of the Einstein condition. A basic question
  is when such a structure is closed, or equivalently locally
  conformally Cotton. In dimension $4$ we obtain a full answer to this
  question and also investigate the incidence of the Bach condition on
  this class of metrics. This is related to earlier results 
  obtained in the Einstein-Weyl context.
\end{abstract}
\maketitle
\tableofcontents

\section{Introduction}
Let $(M^n,g)$ be Riemannian manifold. The metric $g$ is said to be
Einstein if and only if
$$ Ric=\lambda g$$
for some real constant $\lambda$, where $Ric$ is the so-called Ricci
curvature tensor of the metric $g$. Einstein metrics have long had a
priviliged role in geometry. Toward the study of Einstein structures,
and also because Einstein metrics may be obstructed (for example
topologocally), various generalisations of the Einstein condition are
important \cite{besse}. One consists in requiring the so-called {\em
  Cotton tensor} $C$ of the metric $g$ to vanish and weakening further
we might simply require that $g$ be conformal to a Cotton metric.
This is achieved if there is an gradient field
$\zeta$ solving the equation
\begin{equation}\label{pc}
\iota_\zeta W+C=0~,  
\end{equation}
where $\iota_\zeta$ indicates insertion (or interior multiplication) of $\xi$ and $W$ denotes the Weyl tensor of the metric $g$. 

The requirement that $\zeta$ should be exact makes the condition
equation (\ref{pc}) awkward to deal with, directly.  It is not
obvious, for example, how to give local conformal invariants which
characterise metrics which are locally conformally Cotton.  This
suggests that, in the first, instance one might consider a ``Weyl
analogue'' of conformal Cotton equation (\ref{pc}) as follows. We will
say (following \cite{GN}) that a Riemannian manifold is a {\em
  conformal C-space}, if there is a solution to (\ref{pc}), where
$\zeta$ is any section of $TM$. Via a suitable (and natural)
interpretation of $\zeta$ this is a conformally invariant condition. 
In this context we
term (\ref{pc}) the conformal C-space equation.

This move is also motivated by the natural tractor/Cartan structures
of Cartan and Thomas \cite{Cartan,Thomas} (see
\cite{BEGo,CapGotrans,Gau} for modern treatments) and the
corresponding conformal holonomy \cite{Armstrong,Leitner}. Dimension
$n$ conformal Riemannian manifolds are naturally equipped with a rank
$n+2$ vector bundle with a Lorentzian signature metric and compatible
canonical connection. This is the tractor bundle and connection and
the equivalent \cite{CapGotrans} principal bundle structure is the
Cartan connection. An Einstein structure determines a parallel section
$I$ of the standard tractor bundle and, conversely, if a conformal
structure admits a parallel standard tractor field $I$ then this
(parallel section) determines an Einstein metric on an open dense set
(and we say the manifold is {\em almost Einstein} \cite{Goalmost}.)
This is the set where $X(I)$ is non-vanishing, where $X$ is a
canonical homomorphism which takes sections of the tractor bundle to
conformal densities.
 Writing
$\Omega$ for the curvature of the tractor connection, $I$ parallel 
 clearly implies $\Omega I=0$. An obvious weakening
of the almost Einstein condition is to require that there is a section
$I$ of the standard tractor bundle (now not necessarily parallel)
satisfying $\Omega I=0$. 
On the
open set where $X(I)$ is non-vanishing,  $\Omega I=0$ is
exactly (\ref{pc}), the conformal C-space equation
\cite{GN}. Thus in the sense of infinitesimal conformal
holonomy the conformal C-space equation is a vastly weaker requirement
than the Einstein condition (which would have $I$ annihilate the full
jet of the tractor curvature).

Using related ideas it is
straightforward to manufacture conformal obstructions to conformal
C-space metrics.  For example on Riemannian 4-manifolds with
non-vanishing Weyl tensor the conformal invariant
$$
|W|^2C_{abc}+4W^{dijk}C_{ijk}W_{dabc}
$$
(where we have used an obvious abstract index notation) vanishes if
and only if the manifold is a conformal C-space; this result
is easily recovered, or see Proposition 2.5 of \cite{GN}.  

Taking the conformal C-space equation as our basic generalisation of
the Cotton and Einstein conditions, the fundamental question is then
is how far we are from these ``integrable cases''. In the case that
the Weyl curvature is suitably non-degenerate, answers to this may be found in
\cite{GN}, and see also \cite{KTN}. Here we initiate a study of
this issue on Riemannian manifolds, but where the aim is to
remove the assumption of local conditions on the Weyl curvature.  One of
the main results is the following. 
\begin{teo} \label{Main1} If $(M^4,g)$ is a compact conformal
  $C$-space then it is locally conformally Cotton.
\end{teo}

Recall that on 4-manifolds the Bach tensor $B$ (of conformal
relativity \cite{Bach}) is a conformal invariant with leading term
a divergence of the Cotton tensor.  This vanishes on half-flat manifolds, on
locally conformally Einstein structures \cite{besse} and also on the
conformal classes of certain product manifolds \cite{GoLeit}.
Bringing this into the picture leads (see Section 6) to a stronger
result.
\begin{teo} \label{Main2} If $(M^4,g)$ is a compact conformal
  $C$-space which is Bach-flat then it is locally conformally
  Einstein
\end{teo}
\noindent We note that in \cite{KTN} the authors obtained a result that Bach
flat conformal C-spaces are locally conformally Einstein in 4
dimensions, provided that the Weyl tensor satisfies non-degeneracy
conditions. So in the compact Riemannian setting the Theorem improves
their result by removing the need for a non-degeneracy assumption.

Another generalisation of the Einstein condition which has been
studied extensively (e.g.\ \cite{g,pg,ET}) are the Einstein-Weyl
equations. A conformal manifold is said to be Einstein-Weyl if it
admits a compatible torsion-free connection that has vanishing
trace-free symmetrised Ricci curvature. Writing $h$ for the
trace-adjusted (``reduced'') Ricci tensor and $\nabla$ for the
Levi-Civita connection, this problem is equivalent to finding a 1-form
field $\zeta$ and a metric $g$ so that the symmetric part of $h-\nabla
\zeta + \zeta \otimes \zeta$ is pure trace.  There is a close
connection with conformal C-spaces and this plays a role in   the
proofs of the main theorems.

Our paper is organised as follows. In section 2 we collect a number of
basic facts of relevance for the study of conformal $C$-spaces. In
section 3 we study suitablly defined symmetries of an algebraic Weyl
curvature tensor and show how these can be fully understood in
dimension $4$. To prove Theorem \ref{Main1} we first show in section
$4$ that it holds locally, that is on the open subset where the Weyl
tensor does not vanish.  We use a detailed analysis of the properties
of the Weyl curvature of a $4$-dimensional conformal $C$-space
combined with some results from Hermitian geometry \cite{ag}. In
section 5 we establish that the unique continuation property holds for
the class of conformal
C-spaces and this eventually leads to the proof in the compact case. Finally, the last section of the paper is devoted to the proof of Theorem \ref{Main2}.\\

\section{Conformal $C$-spaces and related structures}This section is intended to recall a number of elementary facts concerning the objects we shall subsequently use. Let 
$(M^n,g),n\ge 3$ be a Riemannian manifold, $\nabla$ the Levi-Civita connection associated with the metric $g$ and $R$ the Riemannian curvature tensor, given by 
$R(X,Y)Z=-\nabla_{X,Y}^2Z+\nabla^2_{Y,X}Z$, whenever $X,Y,Z$ are vector fields on $M$. The  Weyl tensor $W$ is defined by the decomposition
\begin{equation}R=W+S
\end{equation}
where the Schouten tensor $S$ is given by $ S=h \bullet g$. Here $h=\frac{1}{n-2}\big(Ric-\frac{s}{2(n-1)}g\big)$  is the reduced Ricci tensor of the metric $g$, whilst the 
Kulkarni-Nomizu product of two symmetric tensors $h$ and $k$ is defined by
$$ (h \bullet k)(x,y,z,t)=h(x,z)k(y,t)+h(y,t)k(x,z)-h(x,t)k(y,z)-h(y,z)k(x,t).$$ The Weyl tensor satisfies the first Bianchi identity
\begin{equation} \label{b1}W(X,Y)Z+W(Y,Z)X+W(Z,X)Y=0.
\end{equation}The second Bianchi identity for $W$ is slightly more complicated and depends on the Cotton tensor $C$, an element of $TM \otimes \Lambda^2(M)$. It is defined by 
$$ C(U,X,Y)=(\nabla_Xh)(Y,U)-(\nabla_Yh)(X,U) $$ 
for all $U,X,Y$ in $TM$ and then 
\begin{equation} \label{b2}
\sigma_{X,Y,Z}\biggl [ (\nabla_XW)(Y,Z,U,T)\biggr ]+(C_U \wedge T-C_T \wedge U)(X,Y,Z)=0
\end{equation}
where $\sigma$ stands for the cyclic sum. An appropriate contraction of the differential Bianchi identity (\ref{b2}) alternatively gives the Cotton tensor $C$ from the Weyl tensor $W$ as
\begin{equation}\delta W= -(n-3)C
\end{equation}
where $\delta W=-\sum \limits_{i=1}^{n}(\nabla_{e_i}W)(e_i,\cdot,
\cdot, \cdot)$ for an arbitrary local orthonormal frame $\{e_i, 1 \le
i \le n\}$ on $M$. The Bach tensor $B$ of the Riemannian manifold
$(M^n,g)$ is defined by
\begin{equation*}
\langle BX,Y \rangle=\sum \limits_{i=1}^{n} \nabla_{e_i}(\delta W)(X,e_i,Y)+W(X,e_i,h(e_i),Y)
\end{equation*}
for all $X,Y$ in $TM$. As it is well known this is symmetric and
tracefree, and moreover in dimension $4$ it is conformally invariant.
In dimension $4$ it also vanishes on (anti)self-dual metrics
\cite{pg}. In this paper we shall mainly study the following class of
Riemannian manifolds.
\begin{defi} \label{def1}Let $(M^n,g)$ be a Riemannian manifold and let $\z$ be a vector field on $M$. Then $(M^n,g,\z)$ is a conformal $C$-space 
if the equation
\begin{equation} \label{Csp}W(\z, \cdot, \cdot, \cdot)+C=0.
\end{equation} is satisfied. If $\z=0$ then $(M^n,g)$ is called a Cotton space.
\end{defi}
It should be noted that conformal $C$-spaces are conformally invariant in the usual sense. Also, a natural sub-class to look at consists in {\it{closed}} 
conformal $C$-spaces, that is conformal $C$-spaces $(M^n,g,\z)$ such that $d\z=0$, which is again a conformally invariant condition. Obviously, closedeness 
in the conformal $C$-space context rephrases globally that a Riemannian metric is locally conformal to that of a Cotton space.
\begin{rema}
Cotton spaces are known to have vanishing Pontriagin classes, fact used in \cite{bou} to show that in $4$ dimensions the non-vanishing of the signature 
implies that the metric is Einstein, provided the manifold is compact. Further results and examples were obtained in dimension $4$ \cite{der} under degeneracy assumptions 
on the spectrum of the Weyl or the Schouten tensor. It is also known that-in the compact case-the metric has to be Einstein when in the presence of a compatible K\"ahler \cite{besse} or closed $G_2$ structure \cite{bryant},\cite{ivanov}. Despite constant interest a complete classification seems to be still missing. 
\end{rema}
Note that in dimensions $n \le 3$ a conformal $C$-space is automatically Cotton due to the absence of algebraic Weyl curvature tensors, hence 
the first interesting dimension in this context is when $n=4$. Related to conformal $C$-spaces are 
Einstein-Weyl structures whose definition we give below.
\begin{defi}\label{EW}
Let $(M^n,g)$ be Riemannian. Then $g$ is Einstein-Weyl  if 
\begin{equation*}
h=\nabla \z-\z \otimes \z-\frac{1}{2}d\z+fg
\end{equation*}
for some vector field $\z$ on $M$ and some smooth function $\z$ on $M$. Moreover $(g,\z)$ is said to be a closed Einstein-Weyl structure if $d\z=0$.
\end{defi}
The question under study in this paper consists in investigating up to what extent a conformal 
$C$-space must be closed. To appproach this, we consider for a given conformal 
$C$-space $(M^n,g,\z)$ the tensor $h_{\z}$ defined by 
\begin{equation*}
h_{\z}=h-\nabla \z+\z \otimes \z
\end{equation*}
and recall that 
\begin{pro} \label{GNl} Let $(M^n,g,\z)$ be a conformal $C$-space. The tensor $h_{\z}$ belongs to the space $\E_W$.
\end{pro}
Therefore one must first understand the algebraic structure of the space $\E_W$ and then explore its 
geometric consequences. In the next section we gather a few general facts to this extent in arbitrary dimensions and 
explore thougroutly the four dimensional case.

\section{Algebraic symmetries of Weyl curvature}
In this section we shall study various algebraic equations akin to produce symmetries of an algebraic Weyl curvature tensor. These are actually insightful when studying various geometric structures, and the relevant 
connections will be made clear in the next section.
\subsection{The various equations}
Let $(V^n,g), n \ge 4$ be  a Euclidean vector space. In what follows shall use the metric to identity without further comment $\otimes^2 V$ with $End(V)$ using the convention $\beta=g(h \cdot, \cdot)$, as well as vectors and 
$1$-forms. As a point of notation, we shall use $\langle \cdot, \cdot \rangle $ for the form inner product induced by $g$. Let $b_1 : \Lambda^2 \otimes \Lambda^2 \to \Lambda^3 \otimes \Lambda^1$ be the Bianchi map 
given by 
\begin{equation*}
(b_1R)(x,y,z)=R(x,y)z+R(y,z)x+R(z,x)y
\end{equation*}
whenever $x,y,z$ belong to $V$ and for any $R$ in $\Lambda^2 \otimes \Lambda^2$, where the standard notation applies.
Consider now a non-vanishing algebraic Weyl-curvature tensor $W$ on $V$ that is an element $W$ of $\Lambda^2 \otimes \Lambda^2$ satisfying the first Bianchi identity, that is $b_1(W)=0$ and which is moreover 
trace-free in the sense that $\sum \limits_{i=1}^n W(\cdot, e_i, \cdot, e_i)=0$ for any orthonormal frame $\{e_i, 1 \le i \le n\}$. Then we can
extend $W$ as a map $W : \otimes^2 V \to \otimes^2 V$ by setting :
$$ W(h)= \sum \limits_{i=1}^n W(e_i, \cdot, he_i, \cdot)$$
for any $h$ in $\otimes^2 V$ and for some arbitrary orthonormal basis $\{e_i\}$ in $V$. This extension of $W$ preserves the tensor type, that is it preserves the splitting $\otimes^2 V=\Lambda^2
\oplus S_0^2 \oplus \mathbb{R}g$. Moreover, the restriction of $W$ to $\Lambda^2(V)$
is given by $<W(v \wedge w), u \wedge q>=W(v,w,u,q)$ for all $v,w,u,q$ in $V$. Using the first Bianchi identity $W$ has to satisfy
this can also be rephrased to say that
$$ W(\alpha)=\frac{1}{2} \sum \limits_{i=1}^{n}W(e_i,Fe_i)$$
for an arbitrary orthonomal basis $\{e_i,1 \le i \le n\}$ and for all $2$-forms
$\alpha$ with associated skew-symmetric endomorphism $F$, that is
$\alpha=g(F \cdot, \cdot)$. We shall be interested in what follows in
the space $\E_W$ of tensors $h$ in $End(V)$ such that
\begin{equation} \label{deg1}
W(x,y,hz,\cdot)+W(y,z,hx,\cdot)+W(z,x,hy, \cdot)=0
\end{equation}
for all $x,y,z$ in $V$. We also define the spaces $\S_W=\E_W \cap S^2, \A_W=\E_W \cap \Lambda^2$ and point out that, a priori, $\E_W$ is not the direct sum of $\S_W$ and $\A_W$. 
The space $\S_W$ has been studied in detail in \cite{bou}. In dimension $4$, as we shall recall later on, additional information is available \cite{der}.
\begin{lema} \label{warmup}
Let $h$ be in $End(V)$. The following hold:
\begin{itemize}
\item[(i)] if $b_1(W(\cdot, \cdot,h\cdot, \cdot))$ belongs to $\Lambda^4$ then $h$ satisfies \eqref{deg1} and 
\begin{equation} \label{deg11}
W(hx,y,z,u)+W(x,hy,z,u)=W(x,y,hz,u)+W(x,y,z,hu)
\end{equation}
whenever $x,y,z,u$ belong to $V$.
\item[(ii)] $h$ satisfies \eqref{deg1} if and only if it satisfies \eqref{deg11}.
\end{itemize}
\end{lema}
\begin{proof}
(i) Let us set $T=b_1(W(\cdot, \cdot,h\cdot, \cdot))$. Then 
\begin{equation*}
W(x,y,hz,u)+W(y,z,hx,u)+W(z,x,hy, u)=T(x,y,z,u)
\end{equation*}
for all $x,y,z,u$ in $V$. We anti-symmetrise in $z,u$ hence 
\begin{equation*}
\begin{split}
&W(x,y,hz,u)+W(x,y,z,hu)+(W(y,z,hx,u)-W(y,u,hx,z)) \\
+&(W(z,x,hy,u)-W(u,x,hy,z))=2T(x,y,z,u)
\end{split}
\end{equation*}
and further 
$$W(x,y,hz,u)+W(x,y,z,hu)-W(hx,y,z,u)-W(x,hy,z,u)=2T(x,y,z,u)$$
after making use of the Bianchi identity. Since $T$ is a four form, it belongs to $S^2(\Lambda^2)$, but since the l.h.s in the equation above belongs to $\Lambda^2(\Lambda^2)$, it must vanish and 
the claim follows.
$\\$
(ii) follows from the Bianchi identity when taking the cylic sum upon $x,y,z$ in \eqref{deg11}.
\end{proof}
\begin{lema} \label{eqv}
The following hold:
\begin{itemize}
\item[(i)] Suppose that $h$ is in $\E_W$. Then
\begin{equation} \label{deg+}
W(hF-Fh^{\star})=W(F)h-h^{\star}W(F)
\end{equation}
and
\begin{equation} \label{deg-}
W(hF+Fh^{\star})=W(F)h+h^{\star}W(F)
\end{equation}
whenever $F$ is a skew symmetric endomorphism of $V$ and where $h^{\star}$ stands for the adjoint of $h$ with respect to the metric $g$.
\item[(ii)] if $h$ in $End(V)$ satisfies both of \eqref{deg+} and \eqref{deg-} then it satisfies \eqref{deg1} as well.
\item[(iii)] the identity \eqref{deg-} is equivalent with \eqref{deg1}.
\end{itemize}
\end{lema}
\begin{proof}
(i) Let us fix an orthonormal basis $\{e_i \}$ in $V$. From (\ref{deg1}) we obtain
\begin{equation*}
W(e_i, Fe_i, hv, w)+W(Fe_i, v, he_i, w)+W(v,e_i,hFe_i,w)=0
\end{equation*}
whenever $v,w$ belong to $V$. After summation, we obtain
\begin{equation*}
2<W(F)hv,w>+\sum \limits_{i=1}^{n} W(Fe_i, v, he_i, w)+W(v,e_i,hFe_i,w)=0
\end{equation*}
Since $\sum \limits_{i=1}^{n} W(Fe_i, v, he_i, w)=-\sum \limits_{i=1}^{n}W(e_i,v,hFe_i,w)=-W(hF)(v,w)$ we find
\begin{equation} \label{mainalg}
W(F)h=W(hF)
\end{equation}
whenever $F$ belongs to $\Lambda^2$. Now the equations in \eqref{deg+}, \eqref{deg-} follow when using that $W$ respects the splitting $\otimes^2V=S^2 \oplus \Lambda^2$.\\
(ii) follows when rewriting \eqref{mainalg} using elements of the form $F=v \wedge w$ where $v,w$ belong to $V$. \\
(iii) again by rewritting \eqref{deg-} by means of decomposable elements of the form $F=z \wedge u$ in $\Lambda^2$ we find 
\begin{equation*}
\begin{split}
&W(z,x,hu,y)-W(u,x,hz,y)+\\
&W(hz,x,u,y)-W(hu,x,z,y)=W(z,u,hx,y)+W(z,u,x,hy)
\end{split}
\end{equation*}
whenever $x,y,z,u$ belong to $V$. The use the Bianchi identity upon the first and third respectively second and fourth terms in the l.h.s. of the equation above shows that $h$ satisfies 
the identity in Lemma \ref{warmup} and therefore it belongs to $\E_W$.
\end{proof}
It remains now to understand up to what extent \eqref{deg+} and \eqref{deg-} are equivalent. 
\begin{lema}\label{eqvlast} Let $W$ be an algebraic Weyl curvature tensor and let $h$ in $\otimes^2 V$ satisfy \eqref{deg+}. Then $h$ equally satisfies \eqref{deg1}. 
\end{lema}
\begin{proof} 
As in the proof of (iii) of the Proposition above we rewrite \eqref{deg+} for decomposable $F$'s in $\Lambda^2$, of the form $F=x \wedge y$, where $x,y$ in $V$. Since 
$$<W(hF)z,u>=W(x,z,hy,u)-W(y,z,hx,u) $$
and 
$$<W(Fh^{\star})z,u>=W(hx,z,y,u)-W(hy,z,x,u) $$
for all $z,u$ in $V$, we arrive at 
\begin{equation} \label{halfdeg}
\begin{split}
&-<W(x,z)u+W(x,u)z,hy>\\
&+<W(y,z)u+W(y,u)z,hx>=W(x,y,hz,u)-W(x,y,z,hu)
\end{split}
\end{equation}
for all $x,y,z,u$ in $V$. Let $T$ in $\Lambda^3 \otimes \Lambda^1$ be defined by $T=b_1(W(\cdot, \cdot,h\cdot, \cdot))$. We rewrite then \eqref{halfdeg} as 
\begin{equation*}
\begin{split}
T(x,y,z,u)=&W(x,y,z,hu)-W(x,u,z,hy)+W(y,u,z,hx)\\
=&-T(x,y,u,z)
\end{split}
\end{equation*}
for all $x,y,z,u$ in $V$. It follows that $T$ belongs to $\Lambda^4$ and we conclude by means of Lemma \ref{warmup}.
\end{proof}
\begin{rema}
It is easy to see that \eqref{halfdeg}, and therefore \eqref{deg1}, is yet equivalent with 
\begin{equation} \label{deg+sym}
W(hS-Sh^{\star})=W(S)h-h^{\star}W(S)
\end{equation}
for all $S$ in $S^2$.
\end{rema}
\begin{lema} \label{ker}
We have $\pi^{-} \E_W \subseteq Ker(W_{\vert \Lambda^2})$, where $\pi^{-} : \otimes^2 V \to \Lambda^2$ is the orthogonal projection.
\end{lema}
\begin{proof} Follows by taking the trace of \eqref{deg1} in the last two arguments.
\end{proof}
Therefore $Ker(W_{\vert \Lambda^2})$ appears as a first obstruction to the equality of  $\E_W$ and $\S_W$ for when $W_{\vert \Lambda^2}$ is injective the latter spaces coincide. We shall show now
that the algebraic structure of $\E_W$ is related to the symmetry group 
$$G_W=\{\gamma \in GL(V): W(\gamma \cdot, \gamma \cdot, \gamma \cdot, \gamma \cdot )=W\}$$ of the Weyl tensor $W$. The Lie algebra $\mathfrak{g}_W$ of 
$G_W$ consists in the space of tensors $h$ in $\otimes^2V$ satisfying 
\begin{equation} \label{lie}
\begin{split}
&W(hx,y,z,u)+W(x,hy,z,u)+\\
&W(x,y,hz,u)+W(x,y,z,hu)=0
\end{split}
\end{equation}
for all $x,y,z,u$ in $V$. Before making explicit the relationship between $\E_W$ and $\mathfrak{g}_W$ we need to establish to reinterpret the identity \eqref{lie} as it has been done for 
\eqref{deg1} and establish the analogous equivalences.
\begin{lema} \label{eqvlie} The following are equivalent:
\begin{itemize}
\item[(i)] $h$ belongs to $\mathfrak{g}_W$
\item[(ii)] $W(hF+h^{\star}F)=-W(F)h-h^{\star}W(F)$ for all $F$ in $\Lambda^2$
\end{itemize}
\end{lema}
\begin{proof} 
Follows by using, with minor changes, the same ingredients as in the proof of Lemma \ref{eqv}. Details are left to the reader.
\end{proof}
\begin{pro} \label{liestr}
We have $[\E_W, \E_W] \subseteq \mathfrak{g}_W$.
\end{pro}
\begin{proof}
We shall make essentially use of the equation \eqref{mainalg}, all tensors in $\E_W$ must satisfy. Let therefore $h_1$ and $h_2$ belong to $\E_W$. Using \eqref{mainalg} we have 
$W(h_1\pi^{-}(h_2F))=W(\pi^{-}(h_2F))h_1$ for all $F$ in $\Lambda^2$. Since $\pi^{-}(h_2F)=\frac{1}{2}(h_2F+Fh_2^{\star})$ we get 
\begin{equation*}
\begin{split}
W(h_1h_2F)+W(h_1Fh_2^{\star})=&W(h_2F+Fh_2^{\star})h_1\\
=&\biggl [ W(F)h_2+h_2^{\star}W(F) \biggl ] h_1
\end{split}
\end{equation*}
after making use of \eqref{mainalg} for $h_2$. Similarly, $W(h_2h_1F)+W(h_2Fh_1^{\star})=\biggl [ W(F)h_1+h_1^{\star}W(F) \biggl ] h_2$ whence 
\begin{equation*}
W([h_1,h_2]F)+W(h_1Fh_2^{\star}-h_2Fh_1^{\star})=-W(F)[h_1,h_2]+h_2^{\star}W(F)h_1-h_1^{\star}W(F)h_2
\end{equation*}
for all $F$ in $\Lambda^2$. But $h_1Fh_2^{\star}-h_2Fh_1^{\star}$ and $h_2^{\star}W(F)h_1-h_1^{\star}W(F)h_2$ are symmetric tensors therefore $\pi^{-}W([h_1,h_2]F)=
-\pi^{-}W(F)[h_1,h_2]$ and this leads to 
\begin{equation*}
W([h_1,h_2]F+F[h_1,h_2]^{\star})=-W(F)[h_1,h_2]-[h_1,h_2]^{\star}W(F)
\end{equation*}
for all $F$ in $\Lambda^2$. Therefore, by using the equivalence of (ii) and (i) in Lemma \ref{eqvlie} we find that $[h_1,h_2]$ belongs to $\mathfrak{g}_W$.
\end{proof}
As it has been done for the space $\E_W$ the obstruction for $\mathfrak{g}_W$ to be contained in $\mathfrak{so}(V)$ is measured as follows.
\begin{lema} \label{ker2}
We have $\pi^{+} \mathfrak{g}_W \subseteq Ker(W_{\vert S^2})$ where $\pi^{+} : \otimes^2V \to S^2$ denotes the orthogonal projection.
\end{lema}
\begin{proof} 
Follows by taking the trace of the identity \eqref{lie} in the variables $x$ and $z$.
\end{proof}
The space $\E_W$ has also an algebraic structure of its own, though different than that of $\mathfrak{g}_W$. Let $\{\cdot, \cdot \} :\otimes^2 V \times \otimes^2 V \to \otimes^2$ denote the anti-commutator. 
\begin{pro} \label{algstr} Let $W$ be an algebraic Weyl curvature tensor. Then $\{\E_W, \E_W \} \subseteq \E_W$.
\end{pro}
\begin{proof}
Given $h$ in $\E_W$ it is enough to show that $h^2$ still belongs to $\E_W$. Now if $F$ belongs to $\Lambda^2$ and since $h$ satisfies \eqref{deg+} we have 
$W(hF-Fh^{\star})=W(F)h-h^{\star}W(F)$ and using this for $\pi^{-}(hF)=\frac{1}{2}(hF+Fh^{\star})$ we obtain 
\begin{equation*} 
\begin{split}
W(h(hF+Fh^{\star})-(hF+Fh^{\star})h^{\star})=&W(hF+Fh^{\star})h-h^{\star}W(hF+Fh^{\star})\\
=& (W(F)h+h^{\star}W(F))h-h^{\star}(W(F)h+h^{\star}W(F))\\
=&W(F)h^2-(h^2)^{\star}W(F)
\end{split}
\end{equation*}
for all $F$ in $\Lambda^2$, where we have used that $h$ satisfies \eqref{deg-}. It follows that $h^2$ satisfies \eqref{deg+} hence the claim follows by making use of Lemma \ref{eqvlast}.\\
\end{proof}
\begin{coro} \label{cor1}
Let $h$ belong to $\E_W$. The following hold:
\begin{itemize}
\item[(i)] $W(hFh^{\star})=h^{\star}W(F)h$ for all $F$ in $\Lambda^2$.
\item[(ii)] $W(hx,hy,z,u)=W(x,y,hz,hu)$ whenever $x,y,z,u$ belong to $V$.
\end{itemize}
\end{coro}
\begin{proof}
(i) By Lemma \ref{eqv} we know that $h$ satisfies \eqref{deg-}, that is $W(hF+Fh^{\star})=W(F)h+h^{\star}W(F)$ for all $F$ in $\Lambda^2$. Using 
this for $\pi^{-}(hF)=\frac{1}{2}(hF+Fh^{\star})$ we compute 
\begin{equation*} 
\begin{split}
W(h(hF+Fh^{\star})+(hF+Fh^{\star})h^{\star})=&W(hF+Fh^{\star})h+h^{\star}W(hF+Fh^{\star})\\
=& (W(F)h+h^{\star}W(F))h+h^{\star}(W(F)h+h^{\star}W(F))\\
=&W(F)h^2+(h^2)^{\star}W(F)+2h^{\star}W(F)h
\end{split}
\end{equation*}
for all $F$ in $\Lambda^2$. Since $h^2$ belongs to $\E_W$ by Proposition \ref{algstr}and therefore satisfies \eqref{deg-} the claim follows.\\
(ii) follows when rewriting (i) by means of decomposable elements of $\Lambda^2$.
\end{proof}
\subsection{The $4$-dimensional case}
Let $(V^4,g)$ be a four dimensional, oriented, Euclidean vector space together with an algebraic Weyl tensor $W$. We consider the splitting  $\Lambda^2(V)=\Lambda^{+} \oplus \Lambda^{-}$ in 
its self-dual resp. anti-self-dual components.
Accordingly, we have the splitting of the Weyl tensor as
$$ W=W^{+}+W^{-}
$$
in its self-dual, resp. anti-self-dual parts. Then $W^{\pm}$ belong to $S^2_{0}(\Lambda^{\pm})$ and let us denote by $\Sigma^{\pm}=\{\lambda_k^{\pm},
1 \le k \le 3\}$ their spectra. Of course $\lambda_1^{\pm}+\lambda_2^{\pm}+\lambda_3^{\pm}=0$. Consider now the corresponding (normalized)
system of eigenforms $W^{\pm}\omega^{\pm}_k=\lambda_k^{\pm} \omega_k^{\pm}, k=1,2,3$. These forms are associated to $g$-compatible
almost complex structures $J_k^{\pm}, 1 \le k \le 3$ that is $\omega_k^{\pm}=g(J_k^{\pm} \cdot, \cdot)$. The almost complex structures
satisfy the quaternion identities i.e. $J_1^{\pm}J_2^{\pm}+J_2^{\pm}J_1^{\pm}=0,
J_3^{\pm}=J_1^{\pm}J_2^{\pm}$ and moreover $[J_k^{+}, J_p^{-}]=0$ for all $1 \le k,p \le 3$.
\begin{lema} \label{shift}Let $h$ in $S^2_0(V)$. The $\{h, J_k^{\pm}\}$ belongs to $\Lambda^{\mp}$ for all $1 \le k \le 3$.
\end{lema}
It will be important for subsequent computations to note that $\vert
\omega_k^{\pm}\vert=2, 1 \le k \le 3$ (here we use the norm on forms).
Define now the endomorphisms
$\sigma_{i,j}=J^{+}_iJ^{-}_j, 1 \le i, j \le 3$; then the
$\sigma_{i,j}$'s are orthogonal involutions of $V$, producing and
orthogonal basis in $S_0^2(V)$. Note also that $\vert \sigma_{i,j} \vert=2, 1 \le i, j \le 3$, where the inner product on $S^2(V)$ is define as usually : $<S_1, S_2>= \sum \limits_{i=1}^{4}
<S_1e_i, S_2e_i>$, for some orthonormal basis $\{e_i, 1 \le i \le 4\}$ in $V$.
\begin{lema} \label{eigsym}We have $W(\sigma_{i,j})=(\lambda_i^{+}+\lambda_j^{-}) \sigma_{i,j}$ for all $1 \le i,j \le 3$.
\end{lema}
\begin{proof}
Let $S$ be in $S^2_0(V)$. Let $\{e_i, 1 \le i \le 4\}$ be an orthonormal basis in $V$ and let $v,w$ be
arbitrary vectors in $V$. We compute by expanding $e^i \wedge v$ in the basis $\omega_k^{\pm}, 1 \le k \le 3$
\begin{equation*}
\begin{split}
W(e_i, v,Se_i,w) &=\frac{1}{2}\sum \limits_{k=1}^{3} \omega_k^{+}(e_i,v)<W(\omega_k^{+})Se_i,w>+
\omega_k^{-}(e_i,v)<W(\omega_k^{-})Se_i,w>\\
&=\frac{1}{2}\sum \limits_{k=1}^{3}\lambda_k^{+} \omega_k^{+}(e_i,v)\omega_k^{+}(Se_i,w)+
\lambda_k^{-}\omega_k^{-}(e_i,v)\omega_k^{-}(Se_i,w).
\end{split}
\end{equation*}
Summing now over $i$ we obtain that
\begin{equation*}
\begin{split}
W(S) &=\sum \limits_{i=1}^{4} W(e_i, \cdot ,Se_i, \cdot ) \\
&=-\frac{1}{2}\sum \limits_{k=1}^{3} \lambda_k^{+}J_k^{+}SJ_k^{+}+\lambda_k^{-}J_k^{-}SJ_k^{-}
\end{split}
\end{equation*}
We now take $S=\sigma_{i,j}, 1 \le i,j \le 3$ to arrive, after a short computation to the proof of the Lemma.
\end{proof}
\begin{coro} \label{supersym}
Any algebraic Weyl tensor is, in $4$-dimensions, subject to the algebraic identities
\begin{equation*} \label{wa}
W\{F,G\}=\{W(F),G\}_{0}+\{F, W(G)\}_{0}
\end{equation*}
and 
\begin{equation*} \label{wa}
-W[F,G]=[W(F),G]+[F, W(G)]
\end{equation*}
whenever $F,G$ are in $\Lambda^2(V)$. 
\end{coro}
\begin{proof}Both claims follow when diagonalising $W$ on $\Lambda^2$ and $S^2$ as mentioned above.
\end{proof}
\begin{pro} \label{split1}The following hold 
\begin{itemize}
\item[(i)] $\A_W=\{\alpha \in \Lambda^2 : W(\alpha)=0\}$
\item[(ii)] $\E_W=\S_W \oplus \A_W$
\end{itemize}
\end{pro}
\begin{proof}
  (i) by Lemma \ref{ker} we need only see that $ker(W_{\vert
    \Lambda^2}) \subseteq \A_W$. Indeed, if $h$ in $\Lambda^2$
  satisfies $W(h)=0$, form the last equation in
  Corollary \ref{supersym} we get $-W[F,h]=[W(F),h]$ for all $F$ in $\Lambda^2$, in other words $h$ satisfies \eqref{deg-}. We conclude now by Lemma \ref{eqv}, (iii).\\
  (ii) Pick $h$ in $\E_W$ and split it as $h=h_s +h_a$ along
  $End(V)=S^2\oplus \Lambda^2$. Then $W(h_a)=0$ by Lemma \ref{ker}
  hence $h_a$ belongs to $\E_W$ by (i) whence so does $h_s$ and the
  proof is finished.
\end{proof}
\begin{rema} When starting from the assumption that $h$ is a Codazzi
  tensor on a Riemannian four manifold a description of the space $\S_W$, together 
  with the constrainst implied on the Weyl tensor
  has been obtained by Derdzinski in \cite{der} by geometric means.
\end{rema}
We finish this section with the following fact, to be used extensively in the next section.
\begin{coro} \label{symker}
Let $h$ be in $\S_W$ such that $Tr(h)$=0.
\begin{itemize}
\item[(i)] if $W^{+}=0$ and $det_{\Lambda^{-}} W^{-} \neq 0$ then we
  must have $h=0$.
\item[(ii)] if $W(h)=0$ and $h$ does not vanish then 
  $Ker(W^{\pm})$ are $1$-dimensional, provided that $W^{\pm} \neq 0$. 
\end{itemize}
\end{coro}
\begin{proof}
(i)It easy to see (see \cite{GN} for instance) that $h$ in $\otimes^2
  V$ belongs to the space $\E_W$ if and only if $W^{\star}(h)=0$.  By
  making use of Lemma \ref{eigsym}, applied to the Weyl curvature
  tensor $W^{\star}$ it follows that $\sum \limits_{1 \le i,j \le
    3}\lambda^{-}_jh_{ij} \sigma_{ij}=0$ where $h=\sum
  \limits_{1 \le i,j \le 3}h_{ij}\sigma_{ij}$ and the claim
  follows eventually.\\
(ii) In this case, from $W^{\star}h=W(h)=0$ one obtains by means of Lemma \ref{eigsym} the system 
\begin{equation*}
\begin{split}
\sum \limits_{1 \le i,j \le 3}(\lambda_i^{+}\mp\lambda^{-}_j)h_{ij} \sigma_{ij}=0
\end{split}
\end{equation*}
where as before $h=\sum \limits_{1 \le i,j \le 3}h_{ij}\sigma_{ij}$. Therefore 
$\sum \limits_{1 \le i,j \le 3}\lambda_i^{+}h_{ij} \sigma_{ij}=\sum \limits_{1 \le i,j \le 3}\lambda^{-}_j h_{ij} \sigma_{ij}=0$. The claim follows 
by taking into account that either $Ker(W^{\pm})$ is $1$-dimensional, situation which clearly leads to $h=0$, or $W^{\pm}=0$.
\end{proof}

\section{A local classification}
We consider in what follows a $4$-dimensional Riemannian manifold $(M^4,g)$ satisfying the conformal C-space condition (\ref{Csp}) for some vector field $\z$.To 
avoid trivial statements we shall assume in what follows that $W$ does not vanish identically, in other words $g$ is not a conformally flat metric. We recall that in $4$-dimensions the well 
known formula holds
\begin{equation} 
\label{4id}
<W_X, W_Y>=\frac{\vert W \vert^2}{4}<X,Y>
\end{equation}for all $X,Y$ in $TM$. Here, $W_X$ stands for the tensor $W(X, \cdot, \cdot, \cdot)$. It 
follows that Weyl nullity vanishes identically in some open set, i.e $W(K, \cdot, \cdot, \cdot)=0$ for some vector field $K$ implies that $K=0$ in the open set of points where $W$ does not 
vanish. Consider now the splitting of $2$-forms
$$\Lambda^2(M)=\Lambda^{-}(M) \oplus \Lambda^{+}(M)$$
in anti-self dual respect. self dual parts. With respect to this 
splitting the Weyl tensor decomposes as $W=W^{-}+W^{+}$. We start by investigating the case when 
\begin{equation} \label{gdeg} 
W(F)=0
\end{equation}
for some two form $F$ on $M$. We split $F=F^{+}+F^{-}$ in its self-dual resp. anti-self-dual components and we consider the open sets 
$D^{\pm}_{F}=\{m \in M : F_m \neq 0\}$ together with $D_F=D^{+}_F \cup D^{-}_F$. Obviously 
\begin{equation} \label{r1}
W^{+}(F^{+})=0
\end{equation}
\begin{equation} \label{r2}
W^{-}(F^{-})=0
\end{equation}
To make statements precise it is also necessary to consider 
$\W^{\pm}=\{m \in M : W^{\pm}_m \neq 0\}$ as well as $\W=\W^{+} \cup \W^{-}=\{m \in M : W_m \neq 0\}$. We work on the open set $D^{+}_F \cap \W^{+}$ which we assume
to be non-empty (actually we will show that this leads to a contradiction so it will turn out that $W^{+}=0$ on $D^{+}_F$ ). Since $W^{+} : \Lambda^{+} \to \Lambda^{+}$ is 
symmetric and trace-free it follows that there are (locally defined, i.e. in some open region around each point in $D^{+}_F \cap \W^{+}$) $g$-compatible almost complex structures $I,J,K$
satisfying the quaternion identities, that is $IJ+JI=0, K=IJ$,  and such that 
\begin{equation} \label{rg1}
W^{+}(\omega_J)=0, \ W^{+}(\omega_I)=\lambda \omega_I \ \mbox{and} \ W^{+}(\omega_K)=-\lambda\omega_K
\end{equation}
for some nowhere vanishing function $\lambda$, locally defined on $D^{+}_F \cap \W^{+}$. Here $\omega_J=g(J\cdot, \cdot)$ etc. in $\Lambda^{+}$ are the so-called K\"ahler forms of the 
almost complex structures above. \\
Let us recall now that the Nijenhuis tensor of the almost 
complex $(g,J)$ is defined by
\begin{equation*}N_J(X,Y)=[X,Y]-[JX,JY]+J[X,JY]+J[JX,Y]
\end{equation*}
whenever $X,Y$ belong to $\Gamma(TM)$. When $N_J=0$ the almost complex $J$ is said to be integrable and actually gives rise to a complex structure. It is customary to call 
$(g,J)$ a Hermitian structure on $M$. It is now a good moment to recall the following important result.
\begin{teo} \label{weylh} \cite{ag} Let $(M^4,g,J)$ be a Hermitian surface. Then the self-dual Weyl tensor $W^{+}$ is given by :
\begin{equation*}W^{+}=\frac{k}{4}(\frac{1}{2}\omega \otimes \omega-\frac{1}{3}Id_{\vert \Lambda^{+}M})+\hat{F}\otimes \omega+\omega \otimes \hat{F}
\end{equation*}
where $k$ is the so-called conformal scalar  curvature and $\hat{F}=\frac{1}{2}d^{+}\theta$. Here $\theta$ is the Lee  form of $(g,J)$, defined 
by $d\omega=\theta \wedge \omega$, where $\omega=g(J\cdot, \cdot)$ is  the so-called K\"ahler form.
\end{teo}
Note that directly from the definition of the Lee form $\theta$ we get after differentation and using that $\omega$ belongs to 
$\Lambda^{+}(M)$ that $\langle \hat{F},\omega \rangle=0$. Theorem \ref{weylh}, giving the structure of the self-dual Weyl tensor of a Hermitian surface, is one of the main ingredients in the proof of the following:
\begin{pro} \label{gdeg2} Let $F$ in $\Lambda^2$ satisfy \eqref{gdeg}. Then $D^{\pm}_{F} \cap \W^{\pm}=\emptyset$, in other words $W^{\pm}$ vanishes identically on 
$D^{\pm}_{F}$.
\end{pro} 
\begin{proof} We will only prove that $D^{+}_{F} \cap \W^{+}=\emptyset$ the proof of the second part of the claim being completely analogous. Therefore, let us assume that 
$D^{+}_{F} \cap \W^{+}$ is not empty and work towards getting a contradiction. The main idea is to show that the almost complex structure $J$ defined in \eqref{rg1} is actually complex.\\
Indeed, it is well known (sse \cite{besse} for instance) that in dimension $4$ integrability (i.e. the vanishing of the Nijenhuis tensor) is equivalent to
\begin{equation} \label{int}(\nabla_{JX}J)JY=(\nabla_XJ)Y
\end{equation}whenever $X,Y$ belong to $TM$, where $\nabla$ denotes the Levi-Civita connection associated with the metric $g$. On the other hand, 
since$(\nabla_XJ)J+J(\nabla_XJ)=0$ for all $X$ in $TM$ we can write 
$$ \nabla J=a_1 \otimes I+a_2 \otimes K $$ 
for some $1$-forms $a_1,a_2$ on $M$. Using that $\omega_J$ anihilates $W^{+}$ we obtain :
$$ \sum \limits_{k=1}^{4} W(e_k, Je_k, X, Y)=0 $$ 
whenever $X,Y$ belong to $TM$ and for some local orthonormal frame $\{e_k, 1 \le k \le 4 \}$. Taking $Y=e_i$ and derivating in the direction of $e_i$ we get
$$ \sum \limits_{1 \le i,k \le 4} (\nabla_{e_i}W)(e_k,Je_k,X,e_i)+W(e_k,(\nabla_{e_i}J)e_k, X, e_i)=0.$$ Since $\delta W=-C$ we get further 
$$ -\sum \limits_{k=1}^4 C(X, e_k, Je_k)+\sum \limits_{k=1}^4 W(\nabla_{e_i}J)(X,e_i)=0.$$
By the conformal C-space equation the first sum equals 
$$-W(\omega_J)(\z,X)=-W^{+}(\omega_J)(\z,X)=0.$$
We compute
\begin{equation*}
\begin{split}
W(\nabla_{e_i}J)(X,e_i)=& a_1(e_i)W(\omega_I)(X, e_i)+a_2(e_i)W(\omega_K)(X, e_i)\\=&\lambda a_1(e_i)<IX,e_i>-\lambda a_2(e_i)<KX,e_i>.
\end{split}
\end{equation*}
Summing over $i$ it follows that $a_1(I \cdot)=a_2(K \cdot)$ or further $a_2=a_1(J \cdot)$, fact which is clearly equivalent with (\ref{int}).
\end{proof}
We shall now relate the result of Proposition \ref{gdeg2} to properties of the form $\z$, part of the defining data of our conformal C-space $(M^4,g,\z)$. This is done by means of the following observation.
\begin{lema} \label{l1}We have $W(d\z)=0$.
\end{lema}
\begin{proof}Follows from the Proposition \ref{GNl} and Lemma \ref{ker}. 
\end{proof}
If the open subsets $D^{\pm}$ of $M$ are defined as $D^{\pm}=D^{\pm}_{d\z}$ we find immediately from Proposition \ref{gdeg2} the following
\begin{coro} \label{alter} 
Let $(M^4,g,\z)$ be a conformal $C$-space. Then $D^{\pm} \cap \W^{\pm}=\emptyset$, that is $W^{\pm}$ vanishes identically on 
$D^{\pm}$.
\end{coro}
What remains to be dealt with is the behaviour of the Weyl tensor $W^{-}$ on $D^{+}$. In other words we shall work on 
$\W^{-} \cap D^{+}$, which, as before, is to be assumed non-empty. Localising further let us 
define $U=\{m \in \W^{-} \cap D^{+} : detW^{-}_m \neq 0\}$. 
\begin{lema} \label{l13} The metric $g$ is Einstein-Weyl on $U$.
\end{lema}
\begin{proof}
By Proposition \ref{GNl} the tracefree part of the tensor $h_{\z}=h-\nabla \z+\z \otimes \z$ belongs to the space $\E_W$. Since in $4$-dimensions the splitting $\E_W=\S_W \oplus \A_W$ holds 
by Proposition \ref{split1}, it follows that the symmetric, tracefree part of $h_{\z}$ belongs to $\S_W$. But on $U$ we have that $W^{+}=0$ since  
$W^{+}$ vanishes on $D^{+}$ by Corollary \ref{alter} and also that $d^{-}\z=0$ since, again by Corollary \ref{alter}, we know that $d^{-} \z$ vanishes on $\W^{-}$. The claim follows now 
by applying Corollary \ref{symker}, (i) to the symmetric, tracefree part of the tensor $h_{\z}$.
\end{proof}
To study the geometry of $(U,g)$ we start from the following Lemma, part of which summarises the information which has already been obtained.
\begin{lema} \label{recap}On the open sub-set $U$ of $M$ the following hold:
\begin{itemize}
\item[(i)] $W^{+}=0$ 
\item[(ii)] $d^{-}\z=0$, that is the $2$-form $F=d\z$ belongs to $\Lambda^{+}$
\item[(iii)]
$\nabla_XF=2\alpha \wedge X-X \lrcorner (F \wedge \z)+3g(\z,X)F$
for all $X$ in $TU$. Here the $1$-form $\alpha$ is given by $\alpha=h(\z,\cdot)+f g(\z, \cdot)-df$ for some 
smooth function $f$ on $U$.
\end{itemize}
\end{lema}
\begin{proof}(i) and (ii) follow directly from Corollary \ref{alter}.\\
(iii) By Lemma \ref{l13} the conformal $C$-space $(U,g)$ is also Einstein-Weyl, that is 
\begin{equation*} \label{e-w}h=\nabla \z-\zeta \otimes \z-\frac{1}{2}d\z+f \cdot g
\end{equation*}
for some smooth function $f$. By differentiation, we get 
\begin{equation*}
\begin{split}
(\nabla_Xh)Y=& \nabla^2_{X,Y}\z-(\nabla_X \z)Y \z-g(\z,Y)\nabla_X\z-\frac{1}{2}Y \lrcorner \nabla_XF+(Xf)g(Y,\cdot)
\end{split}
\end{equation*}
for all $X,Y$ in $TU$. We now skew-symmetrise in $X,Y$ and obtain that 
\begin{equation*}
\begin{split}
(h \bullet g)(\z,Z,X,Y)=&\frac{1}{2}(\nabla_ZF)(X,Y)-F(X,Y)g(\z,Z)\\
-&g(\z,Y)g(\nabla_X\z,Z)+g(\z,X)g(\nabla_Y\z,Z)\\
+&(X.f)g(Y,Z)-(Y.f)g(X,Z)
\end{split}
\end{equation*}
whenever $X,Y,Z$ belong to $TU$, when using that $R=W+h \bullet g$ and taking into account the conformal $C$-space equation. It suffices now to make use of the Einstein-Weyl equation in the 
second and fourth terms in the expansion of $(h \bullet g)(\z,Z,X,Y)$ to obtain, after computing to some extent, the claimed result.
\end{proof}
\begin{lema} \label{intermed} $U$ is the empty set.
\end{lema}
\begin{proof}
  Arguing by contradiction, let us suppose that $U$ is not empty. We
  write $F \wedge \z=\star \beta$ for some $1$-form $\beta$ on $U$ and
  notice that $X \lrcorner (F \wedge \z)=\star (\beta \wedge X)$ for
  all $X$ in $TU$, since $\lrcorner$ and $\wedge$ are dual operators
  with respect to to the Riemannian metric $g$. Therefore (iii) of
  Lemma \ref{recap} becomes
\begin{equation*} \label{twist1}-\nabla_XF=2\alpha \wedge X+\star(X \wedge \beta)+g(\z,X)F
\end{equation*}
for all $X$ in $TU$. Since $F$ belongs to $\Lambda^{+}$ by Lemma \ref{recap}, (ii), it follows that 
\begin{equation} \label{twist2}-\nabla_XF=(X \wedge \gamma)^{+}+g(\z,X)F
\end{equation}whenever 
$X$ belongs to $TU$, where the $1$-form $\gamma$ is given by $\gamma=\frac{1}{2}(2\alpha+\beta)$. From 
the definition of $U$ the $2$-form $F$ is nowhere 
vanishing on $U$. Therefore, we can write $F=\lambda g(J\cdot, \cdot)$ 
on $U$, for some smooth, non-where vanishing, function $\lambda$ and 
$g$-compatible almost complex structure $J$. Localising even further if necessary, we also choose a $g$-compatible almost complex $I$ s.t. 
$IJ+JI=0$ and set $K=IJ$. Then equation (\ref{twist2}) 
becomes  
\begin{equation*}\nabla_X(\lambda J)=\frac{1}{4}\biggl [\gamma(JX)J+\gamma(IX)I+\gamma(KX)K \biggr ]+
\lambda g(\z,X) J
\end{equation*}for all $X$ in $TU$. Identifying the $J$-invariant, resp. $J$-anti-invariant components in the
equation above yields 
\begin{equation} \label{twist22}d\lambda=\frac{1}{4}J\gamma+\lambda \z
\end{equation}and
\begin{equation*} \label{twist3}\lambda \nabla_XJ=\frac{1}{4} \biggl [ \gamma(IX)I+\gamma(KX)K \biggr ]
\end{equation*}
for all $X$ in $TU$. As usually, let us write now $\omega_J, \omega_I,
\omega_K$ for the self-dual $2$-forms associated with $J,I,K$. Also,
we let an endomorphism $G$ of $TU$ act on a $1$-form $\eta$ in
$\Lambda^1(U)$ by $G\eta=\eta(G \cdot)$.  The last equation above
leads to $4\lambda d\omega_J=(I\gamma) \wedge \omega_I+(K\gamma)
\wedge \omega_K$ and further to
$d\omega_J=\frac{\lambda^{-1}}{2}(J\gamma) \wedge \omega_J$, after
making use of the simple algebraic fact that
\begin{equation*}
J\gamma \wedge\omega_J=I \gamma \wedge \omega_I=K \gamma \wedge \omega_K.
\end{equation*} 
But since $F=\lambda\omega_J$ is a closed $2$-form we get $d\omega_J=-\lambda^{-1}d\lambda \wedge \omega_J$ hence $J\gamma=-2d\lambda$. Infering this in 
(\ref{twist22}) gives $\z=\frac{3}{2}\lambda^{-1}d\lambda$, therefore $d\z=0$, so that $F=0$. This contradicts the non-vanishing of $F$ on $U$ hence the proof is complete.
\end{proof}
We are now in position to clarify, locally, up to what extent a four dimensional conformal $C$-space is closed.
\begin{teo} \label{mainlocal}
Let $(M^4,g,\z)$ be a conformal $C$-space. On the open set $\W$ where the Weyl tensor $W$ does not vanish we have that $d\z=0$. In other words, 
the metric $g$ is locally conformal to a Cotton metric on $\W$.
\end{teo}
\begin{proof}
Lemma \ref{intermed} implies that $det(W^{-})=0$ on $\W^{-} \cap D^{+}$. Supposing that $\W^{-} \cap D^{+}$ is not empty it follows that 
$Ker(\W^{-})$ is $1$-dimensional. Proposition \ref{gdeg2} ensures then the vanishing of $W^{-}$, a contradiction. Therefore $\W^{-} \cap D^{+}$ is empty and similarly one shows that $\W^{+} \cap D^{-}$ is empty as well whence the proof of the claim.
\end{proof}

\section{A Weitzenb\"ock formula and the unique continuation of the Weyl tensor}
Let $(M^{n},g,\z), n \ge 3$ be a conformal $C$-space, where $M$ is supposed to be connected. We will produce a 
Weitzenb\"ock type formula for the Weyl tensor which will enable us to prove the unique continuation property 
for $W$ (note that this is well known for Cotton spaces\cite{besse}). Recall that a smoooth section of some vector bundle $E \to M$ has the (strong) unique continuation property if it 
vanishes over $M$ as soon as it vanishes over some non-empty subset of $M$. The weak unique continuation 
property is said to hold if the section vanishes identically as soon as it has, at some point, an 
infinite order contact with the zero section. We shall compute first the action of 
Laplacian on $W$. Note that, as usually, $\Delta=d \delta+\delta d$ when $W$ is considered as a $2$-form with values in $\Lambda^2(M)$.\begin{pro} \label{laplW}For any conformal $C$-space $(M^{n},g,\z), n \ge 3$ where $M$ is supposed to 
be compact the following estimate holds pointwisely
\begin{equation*}\Vert \nabla^{\star} \nabla W \Vert^2 \le C(\Vert \nabla W \Vert^2+\Vert W \Vert^2)
\end{equation*}
for some positive constant $C$.
\end{pro}
\begin{proof}Using the second Bianchi identity (see \ref{b2}) for $W$ we obtain  
\begin{equation*}
\begin{split}(\delta d W)(Y,Z,U,T) &=(\nabla_TC)(U,Y,Z)-(\nabla_UC)(T,Y,Z) \\
&+(\hat{\delta}C)(U,Y) g(Z,T)-(\hat{\delta}C)(T,Y) g(Z,U)\\&-(\hat{\delta}C)(U,Z) g(Y,T)+(\hat{\delta}C)(T,Z) g(Y,U)
\end{split}\end{equation*}
where we have set $\hat{\delta}C=\sum \limits_{i=1}^{n} (\nabla_{e_i}C)(\cdot, e_i, \cdot)$. On the other side, using the conformal $C$-space equation one shows (see \cite{GN}) that
 \begin{equation*}
\hat{\delta}C=-W(\nabla \z+(n-3)\z \otimes \z).
\end{equation*} 
From the above formula combined with the fact that $W(d\z)=0$ (see Lemma \ref{l1}) 
it follows that $\hat{\delta}C$ is a symmetric tensor hence we end up with 
\begin{equation*}
\begin{split}(\delta d W)(Y,Z,U,T) &=(dC)(T,U,Y,Z)-W(\nabla \z +(n-3)\z \otimes \z) \bullet g
\end{split}
\end{equation*}
But $d\delta W=(n-3)dC$ thus 
\begin{equation} \label{wz1}
\begin{split}(\Delta W)(Y,Z,U,T)&=(n-3)(dC)(Y,Z,U,T)-(dC)(U,T,Y,Z) \\
&-W(\nabla \z + (n-3)\z \otimes \z) \bullet g
\end{split}
\end{equation}
Now the classical Weitzenb\"ock formula asserts that $\Delta W=\nabla^{\star} \nabla W+q(W)$ for some endomorphism $q$ (which for our purposes is not 
necessary to make explicit). Since any linear term in the Weyl tensor can be estimated ($M$ is compact) by $const \Vert W \Vert$ we eventually find by 
making use of the triangular inequality under the form $\vert x +y \vert^2 \le 2(\vert x\vert+\vert y \vert^2)$ that 
\begin{equation*}\Vert \nabla^{\star} \nabla W \Vert^2 \le const (\Vert W \Vert^2+\Vert dC \Vert^2).
\end{equation*}On the other hand, by making use of the conformal $C$-space equation we compute 
\begin{equation*}
\begin{split}
(dC)(X,Y,Z,U)&=(\nabla_XC)(Y,Z,U)-(\nabla_YC)(X,Z,U)\\
&=-(\nabla_XW)(\z, Y,Z,U)+(\nabla_YW)(\z,X,Z,U) \\
&-W(\nabla_X{\z},Y,Z,U)+W(\nabla_Y \z,X,Z,U)
\end{split}
\end{equation*}
Making use of the second Bianchi identity in the second line of the above and keeping in mind that $C$ is linear in $W$ 
(because of the conformal $C$-space equation) we find \begin{equation*}dC=-\nabla_{\z}W+\mbox{linear terms in W}.
\end{equation*}Hence 
\begin{equation*}
\begin{split}
\Vert dC \Vert^2 & \le \Vert \nabla_{\z}W\Vert^2+const \Vert W \Vert^2 \le \Vert \z \Vert^2 \Vert \nabla W \Vert^2+const \Vert W 
\Vert^2 \\& \le const (\Vert \nabla W \Vert^2+\Vert W \Vert^2)
\end{split}
\end{equation*}and the result follows immediately.
\end{proof}Using clasical results of Aronszajn \cite{aro,kaz} yields then
\begin{teo} \label{unqW}Let $(M^{n},g,\z)$ be a compact, conformal $C$-space. The Weyl tensor $W$ has the strong unique continuation property.
\end{teo}
$\\$ {\bf{Proof of Theorem \ref{Main1}}}: From Theorem
\ref{mainlocal} we know that $W$ vanishes on the open set $D=\{m
\in M: (d\z)_m \neq 0\}$. If $D$ is empty then trivially $d\z=0$; if not,
by using the unique continuation property established above we obtain
that $W=0$, hence the proof is finished. \quad $\Box$

\section{Bach flat manifolds}
We are interested here in conformal C-spaces $(M^4,g,\z)$ which are
in addition assumed to be Bach flat, that is $B=0$. This has been
previously studied, under some genericity assumptions in \cite{KTN}
and our objective in this section is to discuss the general case of
this setting.
\begin{teo} \label{KNT} Let $(M^4,g,\z)$ be a conformal C-space such
  that $B=0$. The following hold:
\begin{itemize}
\item[(i)] on the open set where $W$ does not vanish $g$ is a 
closed Einstein-Weyl metric.
\item[(ii)] if $M$ is compact and not locally conformally flat then
  $g$ is a closed Einstein-Weyl metric all over $M$.
\end{itemize}
\end{teo} 
\begin{proof} 
  (i) Differentiating the C-space equation and using the definition of
  the Bach tensor one finds \cite{GN,KTN}
that the tensor $h_{\z}$ satisfies the
  additional condition
\begin{equation*}
W(h_{\z})=0. 
\end{equation*}
We also recall that by Theorem \ref{mainlocal}, we have that $d\z=0$
on $\W$, in other words $h_{\z}$ is symmetric on $\W$ or on $M$ if the
latter is compact and in both cases $h_{\z}$ belongs to $\S_{W}$.
Now let $h^0_{\z}$ denote the tracefree part of $h_{\z}$.  Working now on
$\W^{+}$ we see from Corollary \ref{symker}, (ii) that around each
point where $h^0_{\z}$ does not  vanish $W^{+} : \Lambda^{+} \to
\Lambda^{+}$ has $1$-dimensional kernel. But this is a contradiction in view of Proposition \ref{gdeg2}. 
Therefore $h^0_{\z}$ vanishes on $\W^{+}$.  By a similar argument,
$h_{\z}^0$
vanishes on $\W^{-}$ and this proves the claim in (i). \\
(ii) We reinterpret (i) to say that the Weyl tensor $W$
vanishes on the open set where $h^0_{\z}$ does not. If the latter is
assumed non-empty, the use of the unique continuation property for $W$
(cf. Theorem \ref{unqW}) yields the vanishing of $W$. Hence if $g$ is
not locally conformally flat we must have $h^0_{\z}=0$ on $M$, and the
claim follows.
\end{proof}
Under the assumptions above it is easy to see that the metric $g$ must
be locally conformally Einstein. We note that, as shown by Pedersen \&
Swann \cite{PS}, it is sufficient to have Bach flat and the
Einstein-Weyl equations holding to conclude that a compact four
manifold is locally conformally Einstein.
In the compact 4-manifold setting, it is also the case that an
Einstein-Weyl metric with Cotton tensor zero is necessarily locally
conformally Einstein \cite{stefan}. As a variation of this theme we
note the following result.
\begin{teo} \label{KTNglob} Let $(M^4,g,\z)$ be a compact conformal $C$-space. If $g$ is assumed to be Bach flat then either: \\
  (i) $g$ is conformal to an Einstein metric,  \\
  or \\
  (ii) $g$ is a locally conformally flat metric, that is $W=0$.
\end{teo}
\begin{proof} Using the conformal invariance of the conformal C-space
  equation (see \cite{GN} for instance) and that $\z$ is a Weyl one
  form, we shall work in a Gauduchon gauge \cite{g}, that is the
  unique metric in the conformal class of $g$ such that $d^{\star}
  \z=0$. Then a well known result of Gauduchon \cite{pg} states that
  $\z$ is a Killing field on $M$ whence $\z$ is parallel, given that
  it is also closed. That $g$ is an Einstein-Weyl metric is thus
  strengthened to $h=-\z \otimes \z+\lambda g$ for some smooth
  function $\lambda$ on $M$. Again from $\nabla \z=0$ we get
  $C_X=d\lambda \wedge X$ for all $X$ in $TM$, and from the fact that
  the Cotton tensor is trace free it follows that $d\lambda=0$. In
  other words $C=0$ leading to $W(\z, \cdot ,\cdot, \cdot)=0$. Now if
  $\z$ is not zero this yields $W=0$ by (\ref{4id}), and otherwise if
  $\z=0$ then from the formula above $h$ is Einstein.
\end{proof} {\bf{Acknowledgements}}: This research has been supported
through UoA grants and the Marsden Grant no. 02-UOA-108 of the Royal
Society of New Zealand. P-A.N. wishes to thank P.\ Nurowski for a number
of useful discussions.

\end{document}